# Product of parabolic cylinder functions involving Laplace transforms of confluent hypergeometric functions

Ridha Nasri*

*Orange Labs, 38-40 avenue General Leclerc, 92794 Issy-les-Moulineaux, France

In this paper, the product of parabolic cylinder functions $D_\nu(\pm z)D_{\nu+\mu-1}(z)$, with different parameters $\mu$ and $\nu$, are established in terms of Laplace and Fourier transforms of Kummer's confluent hypergeometric functions. The provided integral representations are transformed to easily yield Nicholson-type integral forms and used to derive other series expansions for products of parabolic cylinder functions.

**Keywords:** Parabolic cylinder functions; Weber-Hermite's functions; Confluent hypergeometric functions; Laplace and Fourier transform; Nicholson-type integrals.

*AMS Subject Classification*: 33C15; 44A10; 44A20; 30D10; 30E20.

## 1. Introduction

Parabolic cylinder function, called also Weber-Hermite's function [1], arises in different scientific disciplines and is considered in many branches of mathematics. For instance, in number theory, Helfgott has recently utilized it as a tool in his proof of Goldbach's weak conjecture [2]. He developed a new asymptotic expansion of the parabolic cylinder function using saddle point method [3, 4]. It is also widely considered in physics thanks to its strong relation with the Gaussian normal distribution. Actually, it represents the general moments of the normal distribution; see its integral representations in [5, p. 1028]. Parabolic cylinder function is the solution of the differential equation

$$\frac{d^2y}{dz^2} + \left(\nu + \tfrac{1}{2} - \tfrac{z^2}{4}\right)y = 0 \qquad (1)$$

and is entire function of all $z \in \mathbb{C}$ and $\nu \in \mathbb{C}$.

Following the notation in [5–10], we use $D_\nu(z)$ to denote the parabolic cylinder function. Note that other notations exist in literature such as $U(\nu, z) = D_{-\nu-1/2}(z)$ [2–4, 11–13]. It is known that, for nonpositive integer $\nu$, $D_\nu(z)$ and $D_\nu(-z)$ are two independent solutions of the differential equation (1). Also, $D_\nu(z)$ and $D_{-\nu-1}(iz)$ are two independent solutions for any $\nu \in \mathbb{C}$ [6, p. 117, Eq. (10), (11)].

Parabolic cylinder function $D_\nu(z)$ may be expressed in terms of Kummer's confluent

---

*Corresponding author. Email: ridha.nasri@orange.com

hypergeometric functions $\Phi$ and $\Psi$ by

$$\begin{aligned} D_\nu(z) =& 2^{\frac{\nu}{2}} e^{-\frac{z^2}{4}} \Psi\left(-\tfrac{\nu}{2}, \tfrac{1}{2}; \tfrac{z^2}{2}\right) \\ =& 2^{\frac{\nu}{2}} e^{-\frac{z^2}{4}} \left( \frac{\sqrt{\pi}}{\Gamma(\tfrac{1}{2}-\tfrac{\nu}{2})} \Phi\left(-\tfrac{\nu}{2}, \tfrac{1}{2}; \tfrac{z^2}{2}\right) - \frac{\sqrt{2\pi}z}{\Gamma(-\tfrac{\nu}{2})} \Phi\left(\tfrac{1}{2}-\tfrac{\nu}{2}, \tfrac{3}{2}; \tfrac{z^2}{2}\right) \right) \end{aligned} \quad (2)$$

A number of integral representations and series expansions of parabolic cylinder functions are given in [4–6, 11–13]. Moreover, integral representations of the product of two parabolic cylinder functions have been also investigated and several formulas are provided in [5, 6, 14]. Among them, it is known that (see [5, p. 697] or [6, p. 120])

$$D_{-\nu}\left(ze^{i\frac{\pi}{4}}\right) D_{-\nu}\left(ze^{-i\frac{\pi}{4}}\right) = \frac{\sqrt{\pi}}{\Gamma(\nu)} \int_0^{+\infty} J_{\nu-\frac{1}{2}}\left(\frac{t^2}{2}\right) e^{-zt} dt, \quad \Re\nu > 0, \ \Re z \geq 0 \quad (3)$$

where $J_\nu$ is the Bessel's function of the first kind. Other integral representations of the products of $D_\nu$, involving modified Bessel's function of the second kind, are also found in [5, 6].
Furthermore, Durand in [15] derived a Nicholson-type integral for the product of two parabolic cylinder functions, valid for $\Re\nu > -1$ and for all $z \in \mathbb{C}$.

$$\begin{aligned} &D_\nu(z)^2 + \frac{1}{\sin^2 \pi\nu} \left(\cos \pi\nu D_\nu(z) - D_\nu(-z)\right)^2 \\ &= \frac{2\sqrt{2}\Gamma(\nu+1)}{\pi} \int_0^{+\infty} e^{-(2\nu+1)t + \frac{z^2}{2}\tanh t} \frac{dt}{\sqrt{\sinh 2t}} \end{aligned} \quad (4)$$

Malyshev provided in [7] other useful integral representations for the products $D_\nu(z)D_\nu(\pm z)$, valid for $\Re\nu < 0$. He showed that

$$D_\nu(z)D_\nu(\pm z) = \frac{\sqrt{2} e^{-\frac{z^2}{2}}}{\Gamma(-\nu)} \int_0^{+\infty} e^{(2\nu+1)t \mp z^2(e^{2t} \mp 1)^{-1}} \frac{dt}{\sqrt{\sinh 2t}} \quad (5)$$

It is important to notice that the integral representation (5) of $D_\nu(z)D_\nu(-z)$ provided by Malyshev is equivalent to Durand formula (4) because, using the relation between $D_\nu(z)$ and $D_{-1-\nu}(iz)$ in equation (24), it is easy to show that

$$D_\nu(z)^2 + \frac{1}{\sin^2 \pi\nu} \left(\cos \pi\nu D_\nu(z) - D_\nu(-z)\right)^2 = \frac{2\Gamma(\nu+1)^2}{\pi} D_{-1-\nu}(iz) D_{-1-\nu}(-iz) \quad (6)$$

Therefore, replacing $\nu$ by $-1-\nu$ and $z$ by $iz$ in (4) leads to (5) for $D_\nu(z)D_\nu(-z)$.

Very recently, Glasser provided in [8] an integral representation for the parabolic cylinder function product $D_\nu(x)D_\nu(-y)$ with unrelated arguments $x$ and $y$ satisfying $x > y$. In the same direction, Veestraeten [9] has recently derived other integral representations for $D_\nu(x)D_\nu(y)$ and $D_\nu(x)D_{\nu-1}(y)$ with unrelated real arguments $x$ and $y$ satisfying $x + y \geq 0$. He has also pointed out that, using the recurrence relation of the parabolic cylinder function, it is possible to extend his results to have the same kind of integral representations for $D_\nu(x)D_{\nu-n}(y)$, where $n$ is an integer.

Furthermore, Nicholson-type integral (4) is useful in the study of real zeros of parabolic cylinder functions. In their work [1] about the properties of real zeros of any linear combination of solutions of (1), Elbert and Muldoon established a remarkable identity that can be written, after a straightforward transformation, in terms of parabolic cylinder functions as follows

$$D_{-1-\nu}(iz)\frac{\partial D_{-1-\nu}(-iz)}{\partial \nu} - D_{-1-\nu}(-iz)\frac{\partial D_{-1-\nu}(iz)}{\partial \nu}$$
$$= \frac{\sqrt{2}\pi i}{\Gamma(\nu+1)}\int_0^{+\infty} e^{-(2\nu+1)t+\frac{z^2}{2}\tanh t}\operatorname{erf}\left(z(\tfrac{1}{2}\tanh t)^{1/2}\right)\frac{dt}{\sqrt{\sinh 2t}} \qquad (7)$$

where $\operatorname{erf}(z) = 1 - \sqrt{2/\pi}e^{-z^2/2}D_{-1}\left(z\sqrt{2}\right)$ is the error function.

Identity (7), valid for $\Re\nu > -1$, was the key tool to obtain in [1] an explicit formula for the derivative of a zero of Hermite's function with respect to the parameter. Besides, as can be observed, identities (4), (5) and (7) seem to be special cases of a general formula. To the best of my knowledge, no such general identity for the product of parabolic cylinder functions with different parameters exists in literature.

Primarily motivated by the need to generalize identities (3) to (7), we have in this paper three aims. The first is to establish general formulas covering in particular identity (3). As will be shown later, the new formulas explicitly give the product of parabolic cylinder functions with different parameters and involve the Laplace and Fourier transforms of confluent hypergeometric functions. The second aim of this work is to transform the new found identities to yield general Nicholson-type integrals for the product of parabolic cylinder functions with different parameters. These identities are the generalizations of (4), (5) and (7). In addition, we extend the results from the first and second aims to point out some series expansions of the product of parabolic cylinder functions.

## 2. Integral representations of products of parabolic cylinder functions

In what follows, we establish the Laplace and Fourier transforms of confluent hypergeometric functions of argument $x^2$. We show that these integral transforms are explicit and equal to the product of parabolic cylinder functions with different parameters. Moreover, general Nicholson-type integrals for the parabolic cylinder functions products are deduced.

Throughout the paper, we denote by $\Phi$ and $\Psi$ the Kummer's confluent hypergeometric functions of the first and second kind, respectively.

Let $(\nu)_n = \Gamma(\nu+n)/\Gamma(\nu)$ be the Pochhammer symbol, confluent hypergeometric function $\Phi$ is analytic on C and is defined, as in [5], by the series expansion

$$\Phi(\nu,\mu;z) = \sum_{n=0}^{+\infty}\frac{(\nu)_n}{(\mu)_n}\frac{z^n}{n!} \qquad (8)$$

Equation (8) is equivalent to the following integral representation of $\Phi$

$$\Phi(\nu,\mu;z) = \frac{\Gamma(\mu)}{\Gamma(\nu)\Gamma(\mu-\nu)}\int_0^1 e^{zt}t^{\nu-1}(1-t)^{\mu-\nu-1}dt \ , \ 0 < \Re\nu < \Re\mu \qquad (9)$$

On the other hand, confluent hypergeometric function of the second kind $\Psi$ is defined as in [5, 6] by

$$\Psi(\nu, \mu; z) = \frac{\Gamma(1-\mu)}{\Gamma(\nu+1-\mu)} \Phi(\nu, \mu; z) + \frac{\Gamma(\mu-1)}{\Gamma(\nu)} z^{1-\mu} \Phi(\nu+1-\mu, 2-\mu; z) \quad (10)$$

$$= \frac{1}{\Gamma(\nu)} \int_0^{+\infty} e^{-zt} t^{\nu-1} (1+t)^{\mu-\nu-1} dt, \quad \Re\nu > 0, \ \Re z \geq 0 \quad (11)$$

## 2.1. Laplace and Fourier transforms of confluent hypergeometric functions

In this subsection, we derive some useful identities of the product of parabolic cylinder functions involving Laplace and Fourier transforms of confluent hypergeometric functions.

THEOREM 1 *Let $D_\nu$ and $\Phi$ be defined as in (2) and (8), respectively. Let $a \neq 0$ be a purely imaginary complex number, i.e., $\Re a = 0$. Let $\nu$ and $\mu$ be two complex numbers such that $\Re\mu > 0$. Then, for every complex number $z$ satisfying $\Re z > 0$, we have*

$$\frac{a^{\frac{\nu}{2}-\frac{\mu}{2}}}{(-a)^{\frac{\nu}{2}}} D_{-\nu}\left(\frac{z}{\sqrt{-a}}\right) D_{\nu-\mu}\left(\frac{z}{\sqrt{a}}\right) = \frac{1}{\Gamma(\mu)} \int_0^{+\infty} x^{\mu-1} e^{-zx-\frac{1}{2}ax^2} \Phi\left(\nu, \mu; ax^2\right) dx \quad (12)$$

*Proof.* Let $\nu$ and $\mu$ be two complex numbers such that $0 < \Re\nu < \Re\mu$. It follows from the integral representation (9) of $\Phi$ that we have, for every nonnegative real number $x$ and purely imaginary complex number $a$,

$$\frac{1}{\Gamma(\mu)} x^{\mu-1} e^{-\frac{1}{2}ax^2} \Phi\left(\nu, \mu; ax^2\right) = \frac{1}{\Gamma(\nu)\Gamma(\mu-\nu)} \int_0^1 e^{-ax^2(\frac{1}{2}-t)} (xt)^{\nu-1} (x-xt)^{\mu-\nu-1} d(xt)$$

$$= \frac{1}{\Gamma(\nu)\Gamma(\mu-\nu)} \int_0^x e^{-\frac{a}{2}(x^2-2xt)} t^{\nu-1} (x-t)^{\mu-\nu-1} dt \quad (13)$$

The passage to the last integral is just a matter of change of variable, i.e., $xt$ replaced by $t$.
By writing $x^2 - 2xt = (x-t)^2 - t^2$, it is immediate that

$$\frac{1}{\Gamma(\mu)} x^{\mu-1} e^{-\frac{1}{2}ax^2} \Phi\left(\nu, \mu; ax^2\right) = \frac{1}{\Gamma(\nu)\Gamma(\mu-\nu)} \int_0^x e^{\frac{1}{2}at^2 - \frac{1}{2}a(x-t)^2} t^{\nu-1} (x-t)^{\mu-\nu-1} dt$$

$$= f * g(x) \quad (14)$$

where $f * g$ is the Laplace convolution of functions $f$ and $g$, both defined on nonnegative real number $t$ by respectively:

$$f(t) = \frac{t^{\nu-1}}{\Gamma(\nu)} e^{\frac{at^2}{2}} \text{ and } g(t) = \frac{t^{\mu-\nu-1}}{\Gamma(\mu-\nu)} e^{-\frac{at^2}{2}}.$$

It is known that the Laplace transform of a convolution integral is the product of the Laplace transforms of the two functions [5, p. 1108].

Let $\hat{\Phi}(z)$ be the integral in the right hand side of (12). It then follows from (14) that

$$\hat{\Phi}(z) = \mathcal{L}f(z)\mathcal{L}g(z) \tag{15}$$

where $\mathcal{L}f(z)$ is the Laplace transform of function $f$.

From the integral representation of parabolic cylinder functions in [6, p. 119, Eq.(3)], we have

$$\mathcal{L}f(z) = \int_0^{+\infty} e^{-zt} f(t) dt$$
$$= (-a)^{-\frac{\nu}{2}} e^{-\frac{z^2}{4a}} D_{-\nu}\left(\frac{z}{\sqrt{-a}}\right) \tag{16}$$

and

$$\mathcal{L}g(z) = a^{\frac{\nu}{2}-\frac{\mu}{2}} e^{\frac{z^2}{4a}} D_{\nu-\mu}\left(\frac{z}{\sqrt{a}}\right) \tag{17}$$

Making the product of (16) and (17), we complete the proof. ∎

Note that the restriction made on $\nu$ and $\mu$ at the beginning of the proof is required only for the integral representation of confluent hypergeometric function $\Phi$. There is no need for any restriction on $\nu$ in identity (12), it is valid for any $\nu \in \mathbb{C}$ and $\mu \in \mathbb{C}$ such that $\Re\mu > 0$. The restriction made on $a$ is required because the integral in the right hand side of (12) diverges if $\Re a \neq 0$.

THEOREM 2 *Let $D_\nu$ and $\Psi$ be defined as in (2) and (10), respectively. Let $a \neq 0$ be a complex number with $\Re a \geq 0$. Let $\nu$ and $\mu$ be two complex numbers such that $0 < \Re\mu < 2$. Define $\hat{\Psi}(z)$ be the following integral transform*

$$\hat{\Psi}(z) = \int_0^{+\infty} x^{\mu-1} e^{-zx-\frac{1}{2}ax^2} \Psi\left(\nu,\mu; ax^2\right) dx \tag{18}$$

*Then for every complex number $z$, we have*

$$\frac{\sqrt{2}a^{\frac{\mu}{2}} \sin \pi\mu}{\sqrt{\pi}} \hat{\Psi}(z)$$
$$= e^{\frac{\pi i}{2}\mu} D_{-\nu}\left(\frac{iz}{\sqrt{a}}\right) D_{\mu-\nu-1}\left(-\frac{iz}{\sqrt{a}}\right) + e^{-\frac{\pi i}{2}\mu} D_{-\nu}\left(-\frac{iz}{\sqrt{a}}\right) D_{\mu-\nu-1}\left(\frac{iz}{\sqrt{a}}\right) \tag{19}$$

*and consequently*

$$D_{-\nu}\left(\frac{z}{\sqrt{a}}\right) D_{\mu-\nu-1}\left(-\frac{z}{\sqrt{a}}\right) = \frac{2a^{\frac{\mu}{2}}}{\sqrt{2\pi}} \int_0^{+\infty} \sin\left(zx + \frac{\pi\mu}{2}\right) x^{\mu-1} e^{-\frac{1}{2}ax^2} \Psi\left(\nu,\mu; ax^2\right) dx \tag{20}$$

*Proof.* To prove (19), we assume first that $\Re a = 0$ and $\Re z > 0$.

Starting with the representation (10) of $\Psi$ in terms of $\Phi$ and appealing the relation

$\Gamma(\mu)\Gamma(1-\mu) = \pi \csc \pi\mu$, we have the following computation.

$$\hat{\Psi}(z) = \frac{\pi \csc \pi\mu}{\Gamma(\nu+1-\mu)\Gamma(\mu)} \int_0^{+\infty} x^{\mu-1} e^{-zx-\frac{1}{2}ax^2} \Phi\left(\nu, \mu; ax^2\right) dx$$

$$- \frac{\pi \csc \pi\mu}{\Gamma(\nu)\Gamma(2-\mu)} a^{1-\mu} \int_0^{+\infty} x^{1-\mu} e^{-zx-\frac{1}{2}ax^2} \Phi\left(\nu+1-\mu, 2-\mu; ax^2\right) dx \quad (21)$$

Using the result of theorem 1, it directly follows that

$$\hat{\Psi}(z) = \frac{\pi \csc \pi\mu}{\Gamma(\nu+1-\mu)} \frac{a^{\frac{\nu}{2}-\frac{\mu}{2}}}{(-a)^{\frac{\nu}{2}}} D_{-\nu}\left(\frac{z}{\sqrt{-a}}\right) D_{\nu-\mu}\left(\frac{z}{\sqrt{a}}\right)$$

$$- \frac{\pi \csc \pi\mu}{\Gamma(\nu)} \frac{a^{\frac{\nu+1}{2}-\mu}}{(-a)^{\frac{\nu+1-\mu}{2}}} D_{\mu-1-\nu}\left(\frac{z}{\sqrt{-a}}\right) D_{\nu-1}\left(\frac{z}{\sqrt{a}}\right) \quad (22)$$

The integral (18) of $\hat{\Psi}(z)$ is convergent for every complex $a \neq 0$, such that $\Re a \geq 0$. Then, by analytic continuation, identity (22) holds also for any complex $a \neq 0$ satisfying $\Re a \geq 0$ because the right hand side of (22) is valid for every complex number $a \neq 0$. Moreover, it is worth noting that when $\Re a = 0$, the restriction on $\Re z \geq 0$ is required for the convergence of the integral in the right hand side of (18). It is also clear that the integral in (18) is convergent only for $0 < \Re \mu < 2$.

Now, multiplying both sides by $\frac{a^{\frac{\mu}{2}} \sin \pi\mu}{\pi}$ and using the constraint $\Re a \geq 0$ yield

$$\frac{a^{\frac{\mu}{2}} \sin \pi\mu}{\pi} \hat{\Psi}(z)$$

$$= \frac{e^{\frac{\pi i}{2}\nu}}{\Gamma(\nu+1-\mu)} D_{-\nu}\left(\frac{iz}{\sqrt{a}}\right) D_{\nu-\mu}\left(\frac{z}{\sqrt{a}}\right) - \frac{e^{\frac{\pi i}{2}(\nu+1-\mu)}}{\Gamma(\nu)} D_{\mu-1-\nu}\left(\frac{iz}{\sqrt{a}}\right) D_{\nu-1}\left(\frac{z}{\sqrt{a}}\right) \quad (23)$$

For real parameters an arguments $a$ and $z$, the left hand side of (23) is real so also the right hand side. Then, replacing the right hand side of (23) by its complex conjugate is also valid.

The proof of identity (19) results directly from (23). We have to just take the right hand side of (23) and make some transformations using the following known relation between $D_{-\nu}$ and $D_{\nu-1}$ [6, p. 135].

$$e^{\frac{\pi i}{2}(\nu-1)} D_{-\nu}(iz) + e^{-\frac{\pi i}{2}(\nu-1)} D_{-\nu}(-iz) = \frac{\sqrt{2\pi}}{\Gamma(\nu)} D_{\nu-1}(z) \quad (24)$$

Then, it follows from (23) that

$$\frac{a^{\frac{\mu}{2}} \sin \pi\mu}{\pi} \hat{\Psi}(z) = \frac{e^{\frac{\pi i\nu}{2}}}{\sqrt{2\pi}} D_{-\nu}\left(\frac{iz}{\sqrt{a}}\right) \left(e^{\frac{\pi i}{2}(\nu-\mu)} D_{\mu-\nu-1}\left(\frac{iz}{\sqrt{a}}\right) + e^{-\frac{\pi i}{2}(\nu-\mu)} D_{\mu-\nu-1}\left(-\frac{iz}{\sqrt{a}}\right)\right)$$

$$- \frac{e^{\frac{\pi i}{2}(\nu+1-\mu)}}{\sqrt{2\pi}} D_{\mu-1-\nu}\left(\frac{iz}{\sqrt{a}}\right) \left(e^{\frac{\pi i}{2}(\nu-1)} D_{-\nu}\left(\frac{iz}{\sqrt{a}}\right) + e^{-\frac{\pi i}{2}(\nu-1)} D_{-\nu}\left(-\frac{iz}{\sqrt{a}}\right)\right) \quad (25)$$

Simplifying (25), we arrive at identity (19).

To prove (20) from (19), we recall the integral representation (18) of $\hat{\Psi}(z)$, substitute $z$ with $-iz$, equate from each side even parts on $z$ together and then odd parts also together, we find after simplification the following Fourier Cosine and Sine transforms:

$$\frac{4a^{\frac{\mu}{2}}}{\sqrt{2\pi}} \sin \frac{\pi\mu}{2} \int_0^{+\infty} \cos(zx) x^{\mu-1} e^{-\frac{1}{2}ax^2} \Psi\left(\nu, \mu; ax^2\right) dx$$
$$= D_{-\nu}\left(\tfrac{z}{\sqrt{a}}\right) D_{\mu-1-\nu}\left(-\tfrac{z}{\sqrt{a}}\right) + D_{-\nu}\left(-\tfrac{z}{\sqrt{a}}\right) D_{\mu-1-\nu}\left(\tfrac{z}{\sqrt{a}}\right) \quad (26)$$

and

$$\frac{4a^{\frac{\mu}{2}}}{\sqrt{2\pi}} \cos \frac{\pi\mu}{2} \int_0^{+\infty} \sin(zx) x^{\mu-1} e^{-\frac{1}{2}ax^2} \Psi\left(\nu, \mu; ax^2\right) dx$$
$$= D_{-\nu}\left(\tfrac{z}{\sqrt{a}}\right) D_{\mu-1-\nu}\left(-\tfrac{z}{\sqrt{a}}\right) - D_{-\nu}\left(-\tfrac{z}{\sqrt{a}}\right) D_{\mu-1-\nu}\left(\tfrac{z}{\sqrt{a}}\right) \quad (27)$$

Adding (26) to (27) and using that $\sin\frac{\pi\mu}{2}\cos(zx) + \cos\frac{\pi\mu}{2}\sin(zx) = \sin(zx + \frac{\pi\mu}{2})$, we arrive at (20).
It is important to note that in contrast to (18) and (26) where $\mu$ should satisfy $0 < \Re\mu < 2$, the left hand side of (27) converges for $-1 < \Re\mu < 3$. ∎

Recall form [5, p. 1028] that Bessel's functions are special cases of confluent hypergeometric functions.

$$J_{\nu-1/2}(\tfrac{x^2}{2}) = \frac{2^{1-2\nu}}{\Gamma(\nu+\frac{1}{2})} x^{2\nu-1} e^{-\frac{1}{2}ix^2} \Phi\left(\nu, 2\nu; ix^2\right) \quad (28)$$

$$K_{\nu-1/2}(\tfrac{x^2}{2}) = \sqrt{\pi} x^{2\nu-1} e^{-\frac{1}{2}x^2} \Psi\left(\nu, 2\nu; x^2\right) \quad (29)$$

Set then $\mu = 2\nu$, $a = i$ in (12) and $a = 1$ in (20) and use $\Gamma$ duplication formula

$$\Gamma(2\nu) = \frac{2^{2\nu-1}}{\sqrt{\pi}} \Gamma(\nu)\Gamma(\nu+\tfrac{1}{2}) \quad (30)$$

we have immediately the following identities, already known in literature [5, 6].

$$D_{-\nu}\left(ze^{\frac{\pi i}{4}}\right) D_{-\nu}\left(ze^{-\frac{\pi i}{4}}\right) = \frac{\sqrt{\pi}}{\Gamma(\nu)} \int_0^{+\infty} e^{-zx} J_{\nu-\frac{1}{2}}\left(\tfrac{x^2}{2}\right) dx, \quad \Re\nu > 0, \Re z \geq 0 \quad (31)$$

and

$$D_{-\nu}(z) D_{\nu-1}(-z) = \frac{\sqrt{2}}{\pi} \int_0^{+\infty} \sin(zx + \pi\nu) K_{\nu-\frac{1}{2}}\left(\tfrac{x^2}{2}\right) dx, \quad 0 < \Re\nu < 1 \quad (32)$$

*Remark* 1 The integral transform $\hat{\Psi}(z)$ in (18) can be analytically continued for all $\Re\mu > 0$ by just removing the singularity of $\Psi(\nu,\mu;x^2)$ at $x=0$ as follows (for $a=1$).

$$\hat{\Psi}(z) = \int_0^{+\infty} e^{-zx-\frac{x^2}{2}} \left( x^{\mu-1}\Psi\left(\nu,\mu;x^2\right) - \frac{\Gamma(\mu-1)}{\Gamma(\nu)} \sum_{k=0}^n \frac{(\nu+1-\mu)_k}{(2-\mu)_k} \frac{x^{2k+1-\mu}}{k!} \right) dx$$

$$- \frac{\pi \csc \pi\mu}{\Gamma(\nu)} e^{\frac{z^2}{4}} \sum_{k=0}^n \frac{(\nu+1-\mu)_k(2-\mu)_{2k}}{(2-\mu)_k k!} D_{\mu-2-2k}(z) \qquad (33)$$

where $n$ is any natural integer that should satisfy $n \geq \lfloor \frac{\mu-1}{2} \rfloor$.

THEOREM 3 *Let $D_\nu$, $\Phi$ and $\Psi$ be defined as in (2), (8) and (10), respectively. Let $\nu$ be a complex number such that $\Re\nu > 0$. Then, the following integral transforms hold for every complex number $z$.*

$$e^{\frac{z^2}{2}} D_{-\nu}(z)^2 = \frac{2^{2\nu}}{\Gamma(2\nu)} \int_0^{+\infty} x^{2\nu-1} e^{-2zx-2x^2} \Phi\left(\nu, \nu+\tfrac{1}{2}; x^2\right) dx \qquad (34)$$

$$e^{\frac{z^2}{2}} D_{-\nu}(z) D_{-\nu}(-z) = \frac{2}{\Gamma(\nu)} \int_0^{+\infty} \cosh(2zx) x^{2\nu-1} e^{-2x^2} \Psi\left(\nu, \nu+\tfrac{1}{2}; x^2\right) dx \qquad (35)$$

$$e^{-\frac{z^2}{2}} D_{-\nu}(z) D_{-\nu}(-z) = \frac{2}{\Gamma(\nu+\tfrac{1}{2})} \int_0^{+\infty} \cos(2zx) e^{-2x^2} \Phi\left(\tfrac{1}{2}, \nu+\tfrac{1}{2}; x^2\right) dx \qquad (36)$$

*Proof.* To prove (34), let $\Re\nu > 0$ and observe that the function $x \mapsto x^{2\nu-1} e^{-2x^2} \Phi\left(\nu, \nu+\tfrac{1}{2}; x^2\right)$ can be seen as a convolution integral. We have actually form the integral representation (9) of $\Phi$ that

$$\frac{\Gamma(\nu)\Gamma(\tfrac{1}{2})}{\Gamma(\nu+\tfrac{1}{2})} x^{2\nu-1} e^{-2x^2} \Phi\left(\nu, \nu+\tfrac{1}{2}; x^2\right) = x^{2\nu-1} \int_0^1 e^{-x^2(2-t)} \frac{t^{\nu-1}}{\sqrt{1-t}} dt \qquad (37)$$

Making the change of variable $u: t \mapsto x\sqrt{1-t}$ yields

$$\begin{aligned}
\frac{\Gamma(\nu)\Gamma(\tfrac{1}{2})}{\Gamma(\nu+\tfrac{1}{2})} x^{2\nu-1} e^{-2x^2} \Phi\left(\nu, \nu+\tfrac{1}{2}; x^2\right) &= 2\int_0^x e^{-x^2-u^2}(x^2-u^2)^{\nu-1} du \\
&= \int_{-x}^x e^{-\frac{1}{2}(x+u)^2 - \frac{1}{2}(x-u)^2}(x+u)^{\nu-1}(x-u)^{\nu-1} du \\
&= \int_0^{2x} e^{-\frac{u^2}{2} - \frac{1}{2}(2x-u)^2} u^{\nu-1}(2x-u)^{\nu-1} du \\
&= f * f(2x) \qquad (38)
\end{aligned}$$

where $f(t) = t^{\nu-1} e^{-\frac{t^2}{2}}$.

Applying now Laplace transform to the function $x \mapsto [f * f](2x)$ and appealing the integral representation of parabolic cylinder function from [6, p. 119, Eq.(3)], we obtain

$$\frac{2}{\Gamma(\nu)^2} \int_0^{+\infty} e^{-2zx} [f * f](2x) dx = e^{\frac{z^2}{2}} D_{-\nu}(z)^2 \qquad (39)$$

Using the duplication formula (30) of $\Gamma$ function, we complete the proof of identity (34).

To prove (35), we use the following identity:
Let $g$ and $h$ be two summable functions on $[0, +\infty)$, then we have

$$\int_0^{+\infty} g(x)dx \int_0^{+\infty} h(x)dx = \int_0^{+\infty} \int_0^{+\infty} g(x)h(x+t)dxdt + \int_0^{+\infty} \int_0^{+\infty} h(x)g(x+t)dxdt \qquad (40)$$

Its proof is simple. Consider the primitive integral functions $G$ and $H$ of respectively $g$ and $h$:

$$G(x) = \int_x^{+\infty} g(t)dt, \text{ and } H(x) = \int_x^{+\infty} h(t)dt$$

It is clear that $G(+\infty) = H(+\infty) = 0$. Also, observe that the right hand side of (40) can be seen as

$$-\int_0^{+\infty} \left( G'(x)H(x) + G(x)H'(x) \right) dx = G(0)H(0)$$

Applying now relation (40) to $g(t) = t^{\nu-1} e^{-zt - \frac{t^2}{2}}$ and $h(t) = t^{\nu-1} e^{zt - \frac{t^2}{2}}$ yields

$$\Gamma(\nu)^2 e^{\frac{z^2}{2}} D_{-\nu}(z) D_{-\nu}(-z) = 2 \int_0^{+\infty} \cosh(zt) dt \int_0^{+\infty} e^{-\frac{x^2}{2} - \frac{(x+t)^2}{2}} x^{\nu-1}(x+t)^{\nu-1} dx \qquad (41)$$

Changing the variable $x$ with $\frac{t}{2}(\sqrt{u+1} - 1)$ in the inside integral of (41), we have $dx = \frac{tdu}{4\sqrt{u+1}}$ and thus

$$\int_0^{+\infty} e^{-\frac{x^2}{2} - \frac{(x+t)^2}{2}} x^{\nu-1}(x+t)^{\nu-1} dx = \frac{1}{2} \left(\frac{t}{2}\right)^{2\nu-1} \int_0^{+\infty} e^{-\frac{t^2}{4}(u+2)} \frac{u^{\nu-1}}{\sqrt{u+1}} du$$

$$= \frac{\Gamma(\nu)}{2} \left(\frac{t}{2}\right)^{2\nu-1} e^{-\frac{t^2}{2}} \Psi\left(\nu, \nu + \frac{1}{2}; \frac{t^2}{4}\right) \qquad (42)$$

Changing the variable $t$ with $2t$ in (41) and using (42), we complete the proof of (35).

Identity (36) is a direct consequence of (34) and (35). To prove it, we shall use the relation between $\Psi$ and $\Phi$ in (10) and apply some transformations utilizing identity (24). We do not need any restriction on $\nu$ for identity (36); it is valid for every $\nu \in \mathbb{C}$. ∎

I tried to generalize (34) to (36) for the product of parabolic cylinder functions with different parameters $\nu$ and $\mu$ but it seems to do not involve integral transform of confluent hypergeometric functions.

## 2.2. *Nicholson-type integrals for products of parabolic cylinder functions*

As explained in the introduction of this paper, we interest to generalize identities (4) to (7). The idea is to transform the results from the previous theorems to obtain the following corollaries.

COROLLARY 1 *Let $D_\nu$ be defined as in (2). Let $\mu$ and $\nu$ be two complex numbers such that $\Re\nu < 0$. Then we have for $|arg(z)| < \frac{\pi}{4}$*

$$D_\nu(z)D_{\nu+\mu-1}(z) = \frac{2^{1-\frac{\mu}{2}}}{\Gamma(-\nu)} \int_0^{+\infty} e^{(2\nu+\mu)t - \frac{z^2}{4}\coth t} D_{\mu-1}\left(z\sqrt{\coth t}\right) \frac{dt}{(\sinh 2t)^{\frac{\mu}{2}}} \qquad (43)$$

*Proof.* Here, we employ again the result of theorem 1 and we recall the integral representation of confluent hypergeometric function in (9). Assume temporarily that $0 < \nu < \mu$, $a \in i\mathbb{R}^*$ and $z > 0$, we have thus the following direct computation.

$$\frac{a^{\frac{\nu}{2}-\frac{\mu}{2}}}{(-a)^{\frac{\nu}{2}}} D_{-\nu}\left(\frac{z}{\sqrt{-a}}\right) D_{\nu-\mu}\left(\frac{z}{\sqrt{a}}\right)$$

$$= \frac{1}{\Gamma(\nu)\Gamma(\mu-\nu)} \int_0^{+\infty} \int_0^1 x^{\mu-1} e^{-zx - \frac{1}{2}ax^2(1-2t)} t^{\nu-1}(1-t)^{\mu-\nu-1} dt dx$$

$$= \frac{\Gamma(\mu)}{\Gamma(\nu)\Gamma(\mu-\nu)} \int_0^1 e^{\frac{z^2}{4a(1-2t)}} D_{-\mu}\left(\frac{z}{\sqrt{a(1-2t)}}\right) \frac{t^{\nu-1}(1-t)^{\mu-\nu-1}}{(a(1-2t))^{\frac{\mu}{2}}} dt \qquad (44)$$

The last equation is obtained by firstly interchanging the order of integration in the double integral because both converge for the condition made on $\nu$, $\mu$, $a$ and $z$; and secondly by evaluating the inside one using the integral representation of parabolic cylinder function from [6, p. 119, Eq.(3)].
We observe also that the last integral in the right hand side of (44) is convergent for $a$ not necessary purely imaginary complex number assuming of course some restrictions on $z$. Hence, by analytic continuation, the right hand side of (44) is also equal to the left hand side for positive real number $a$. This enables us to set $a = 1$ in what it follows for the proof.
Making the change of variable $x = 1 - 2t$ and splitting the integral into two parts, we have

$$e^{\frac{\pi i \nu}{2}} D_{-\nu}(iz) D_{\nu-\mu}(z) = \frac{\Gamma(\mu) 2^{1-\mu}}{\Gamma(\nu)\Gamma(\mu-\nu)} \int_{-1}^1 e^{\frac{z^2}{4x}} D_{-\mu}\left(\frac{z}{\sqrt{x}}\right) \frac{(1-x)^{\nu-1}(1+x)^{\mu-\nu-1}}{x^{\frac{\mu}{2}}} dx$$

$$= \frac{\Gamma(\mu) 2^{1-\mu}}{\Gamma(\nu)\Gamma(\mu-\nu)} \int_0^1 e^{\frac{z^2}{4x}} D_{-\mu}\left(\frac{z}{\sqrt{x}}\right) \frac{(1-x)^{\nu-1}(1+x)^{\mu-\nu-1}}{x^{\frac{\mu}{2}}} dx$$

$$+ \frac{\Gamma(\mu) 2^{1-\mu}}{\Gamma(\nu)\Gamma(\mu-\nu)} e^{\frac{\pi i \mu}{2}} \int_0^1 e^{-\frac{z^2}{4x}} D_{-\mu}\left(\frac{iz}{\sqrt{x}}\right) \frac{(1+x)^{\nu-1}(1-x)^{\mu-\nu-1}}{x^{\frac{\mu}{2}}} dx \qquad (45)$$

Making once again the change of variable $x = \tanh t$ in both integrals of the right hand side of (45) and simplifying, we arrive at

$$e^{\frac{\pi i \nu}{2}} D_{-\nu}(iz) D_{\nu-\mu}(z) = \frac{2^{1-\frac{\mu}{2}} \Gamma(\mu)}{\Gamma(\nu)\Gamma(\mu-\nu)} \int_0^{+\infty} e^{-(2\nu-\mu)t + \frac{z^2}{4}\coth t} D_{-\mu}\left(z\sqrt{\coth t}\right) \frac{dt}{(\sinh 2t)^{\frac{\mu}{2}}}$$

$$+ \frac{2^{1-\frac{\mu}{2}} \Gamma(\mu)}{\Gamma(\nu)\Gamma(\mu-\nu)} e^{\frac{\pi i \mu}{2}} \int_0^{+\infty} e^{(2\nu-\mu)t - \frac{z^2}{4}\coth t} D_{-\mu}\left(iz\sqrt{\coth t}\right) \frac{dt}{(\sinh 2t)^{\frac{\mu}{2}}} \quad (46)$$

We take imaginary parts from each side of (46) using identity (24). It follows that

$$D_{\nu-1}(z) D_{\nu-\mu}(z) = \frac{2^{1-\frac{\mu}{2}}}{\Gamma(\mu-\nu)} \int_0^{+\infty} e^{(2\nu-\mu)t - \frac{z^2}{4}\coth t} D_{\mu-1}\left(z\sqrt{\coth t}\right) \frac{dt}{(\sinh 2t)^{\frac{\mu}{2}}} \quad (47)$$

Replacing $\nu$ by $\nu + \mu$ gives the desired identity (43).
The integral of (43) is convergent at $t = +\infty$ for $\Re\nu < 0$, while at $t = 0$, it converges provided that $|arg(z)| < \frac{\pi}{4}$ or $\Re\mu < 2$ and $|arg(z)| = \frac{\pi}{4}$. ∎

Equation (43) is the generalization of identity (5) because with $\mu = 1$ in (43), we readily find the Nicholson-type integral of $D_\nu(z)^2$ which was derived by Malyshev in [7, p. 7, Eq. (23)] and extended by Glasser in [8] or by Veestraeten in [9] for unrelated arguments.

In addition, we extend results of theorem 2 to establish the following corollary that generalizes identities (4) and (7).

COROLLARY 2 *Let $D_\nu$ be defined as in (2). Let $\mu$ and $\nu$ be two complex numbers such that $\Re\mu < 2$ and $\Re\nu < 0$. Then, for every complex number $z$, we have*

$$D_\nu(-z) D_{\nu+\mu-1}(z) = \frac{2^{1-\frac{\mu}{2}}}{\Gamma(-\nu)} \int_0^{+\infty} e^{(2\nu+\mu)t - \frac{z^2}{4}\tanh t} D_{\mu-1}\left(z\sqrt{\tanh t}\right) \frac{dt}{(\sinh 2t)^{\frac{\mu}{2}}} \quad (48)$$

*Proof.* We proceed as in (43). To prove (48), we use the result of theorem 2 with $a = 1$ and we recall the integral representation of confluent hypergeometric function $\Psi$ in (11). It follows then from (20) that

$$D_{-\nu}(z) D_{\mu-\nu-1}(-z)$$
$$= \frac{2}{\sqrt{2\pi}\Gamma(\nu)} \int_0^{+\infty} \int_0^{+\infty} \sin\left(zx + \frac{\pi\mu}{2}\right) x^{\mu-1} e^{-\frac{1}{2}x^2(1+2t)} t^{\nu-1} (1+t)^{\mu-\nu-1} dt dx \quad (49)$$

Interchanging the order of integration in (49), evaluating the inside integral (the one with variable $x$) using the integral representation of $D_{\mu-1}(-z)$ in [6, p. 120, Eq.(4)], then making the change of variable $y = \frac{1}{1+2t}$ and simplifying, we obtain

$$D_{-\nu}(z) D_{\mu-\nu-1}(-z) = \frac{2^{1-\mu}}{\Gamma(\nu)} \int_0^1 e^{-\frac{z^2 y}{4}} D_{\mu-1}\left(-z\sqrt{y}\right) \frac{(1-y)^{\nu-1}(1+y)^{\mu-\nu-1}}{y^{\frac{\mu}{2}}} dy \quad (50)$$

Changing the variable of the integration $y = \tanh t$ as in (45) and simplifying, we have

$$D_{-\nu}(z)D_{\mu-\nu-1}(-z) = \frac{2^{1-\frac{\mu}{2}}}{\Gamma(\nu)}\int_0^{+\infty} e^{-(2\nu-\mu)t - \frac{z^2}{4}\tanh t} D_{\mu-1}\left(-z\sqrt{\tanh t}\right)\frac{dt}{(\sinh 2t)^{\frac{\mu}{2}}} \quad (51)$$

Replacing $\nu$ by $-\nu$ and $z$ by $-z$, we complete the proof.
The integral of (48) is convergent at $t = +\infty$ for $\Re\nu < 0$, while at $t = 0$, it converges provided that $\Re\mu < 2$. ∎

*Remark* 2 Using the differential equation (1) satisfied by the parabolic cylinder function and making the first four derivatives of the function $y : z \mapsto D_\nu(z)D_{\nu+\mu-1}(z)$ with respect to $z$, it is easy to prove that the function $y$ satisfies the following fourth order differential equation

$$\frac{d^4 y}{dz^4} + 4\left(\nu + \frac{\mu}{2} - \frac{z^2}{4}\right)\frac{d^2 y}{dz^2} - 3z\frac{dy}{dz} + \mu(\mu-2)y = 0 \quad (52)$$

An alternative proof of corollary 1 is then to show that the right hand side of (43) respects also the differential equation (52) and both sides of (43) are equal for the initial conditions, i.e., $y^{(n)}(0)$, where $0 \le n \le 3$.
Similarly, one can prove that both sides of equation (48) satisfy also the fourth order differential equation (52) and are equal at the initial conditions.

Now, setting $\mu = 1$ in (48) directly provides Nicholson-type integral (5) of $D_\nu(z)D_\nu(-z)$ that was established by Malyshev in [7, p. 7, Eq. (26)] and by Durand in [15] with another form as explained in the introduction.
It would be nice to outline other particular cases of identity (48). Then, if we take from (48) the even and odd parts separately using the expression of $D_\nu$ in (2), we establish the following Nicholson-type integrals:

$$D_\nu(z)D_{\nu+\mu-1}(-z) + D_\nu(-z)D_{\nu+\mu-1}(z)$$
$$= \frac{2\sqrt{2\pi}}{\Gamma(-\nu)\Gamma(1-\frac{\mu}{2})}\int_0^{+\infty}\frac{e^{(2\nu+\mu)t}}{(\sinh 2t)^{\frac{\mu}{2}}}\Phi\left(\frac{\mu}{2}, \frac{1}{2}; -\frac{z^2}{2}\tanh t\right) dt \quad (53)$$

$$D_\nu(z)D_{\nu+\mu-1}(-z) - D_\nu(-z)D_{\nu+\mu-1}(z)$$
$$= \frac{4z\sqrt{\pi}}{\Gamma(-\nu)\Gamma(\frac{1}{2}-\frac{\mu}{2})}\int_0^{+\infty}\frac{e^{(2\nu+\mu)t}}{(\sinh 2t)^{\frac{\mu}{2}}}\sqrt{\tanh}\,\Phi\left(\frac{1}{2}+\frac{\mu}{2}, \frac{3}{2}; -\frac{z^2}{2}\tanh t\right) dt \quad (54)$$

Equation (53) is valid for $\Re\mu < 2$ and $\Re\nu < 0$, while (54) is valid provided that $\Re\mu < 3$ and $\Re\nu < 0$. However, it is possible to remove the restriction on $\mu$ for (53) and (54).
We take the derivative of both sides of (54) with respect to $\mu$ and then set $\mu = 1$, we find

$$D_\nu(z)\frac{\partial D_\nu(-z)}{\partial \nu} - D_\nu(-z)\frac{\partial D_\nu(z)}{\partial \nu} = -\frac{\sqrt{2\pi}z}{\Gamma(-\nu)}\int_0^{+\infty}\frac{e^{(2\nu+1)t}}{\cosh t}\Phi\left(1, \frac{3}{2}; -\frac{z^2}{2}\tanh t\right) dt \quad (55)$$

where $\frac{-2z}{\sqrt{\pi}}\Phi\left(1, \frac{3}{2}; -z^2\right)$ can be reduced to $ie^{-z^2}\text{erf}(iz)$; with erf($z$) is recalled the error function.
Identity (55), valid for $\Re\nu < 0$, is equivalent to [1, p. 8, Eq.(5.13)] provided by Elbert and Muldoon.
Similarly, the derivative of (53) with respect to $\mu$ immediately leads, for $\mu = 2$ after removing singularity, to

$$D_\nu(z)\frac{\partial D_{\nu+1}(-z)}{\partial \nu} + D_\nu(-z)\frac{\partial D_{\nu+1}(z)}{\partial \nu} = \frac{\sqrt{2\pi}z^2}{2\Gamma(-\nu)}\int_0^{+\infty}\frac{e^{2(\nu+1)t}}{\cosh(t)^2}\Phi\left(1,\tfrac{3}{2};-\tfrac{z^2}{2}\tanh t\right)dt$$
$$+ \sqrt{2\pi}\frac{\ln 2 + \psi(-\nu/2)}{2\Gamma(-\nu)} \quad (56)$$

where $\psi(\nu) = \frac{d}{d\nu}\ln\Gamma(\nu)$ is the Digamma function.

### 3. Series expansions of products of parabolic cylinder functions

We show in this section that, using the integral representation from the previous section, we can establish some series expansions of products of parabolic cylinder functions with different parameters.

We use here Tricomi series expansion of confluent hypergeometric function $\Psi$ in terms of Laguerre polynomial [16, p. 275, Eq.(3)].

$$\Psi(\nu,\mu;x) = \frac{1}{\Gamma(\nu)}\sum_{n=0}^{+\infty}\frac{L_n^{\mu-1}(x)}{n+\nu} \quad (57)$$

Identity (57) is valid for positive real $x$ and for $\mu < \frac{3}{2}$. It is also known that Laguerre polynomial is a particular case of confluent hypergeometric function $\Psi$; see [16, p. 268, Eq.(36)].

$$L_n^{\mu-1}(x) = \frac{(-1)^n}{n!}\Psi(-n,\mu;x) \quad (58)$$

Appealing now identity (20) for $a = 1$, representing then $\Psi(\nu,\mu;x^2)$ in terms of Laguerre polynomial as in (57) and finally integrating term by term in the sum, we establish the following series expansion of the product of parabolic cylinder functions.

$$D_{-\nu}(-z)D_{\mu-\nu-1}(z) = \frac{1}{\Gamma(\nu)}\sum_{n=0}^{+\infty}\frac{D_n(z)D_{\mu+n-1}(z)}{n!(n+\nu)} \quad (59)$$

where the series of (59) converges for every $z \in \mathrm{R}$ and for $\mu < \frac{3}{2}$.
Transforming (59) with identity (24), it is easy after simplification to find also

$$D_{-\nu}(z)D_{\mu-\nu-1}(z) = \frac{1}{\Gamma(\nu)}\sum_{n=0}^{+\infty}\frac{(-1)^n D_n(z)D_{\mu+n-1}(z)}{n!(n+\nu)} \quad (60)$$

which is valid for $z \geq 0$ and $\mu < \frac{3}{2}$.
If we let $\mu = 1$ in (59) and (60), we deduce the series expansion of $D_{-\nu}(z)D_{-\nu}(\pm z)$ in terms of Hermite's polynomials[1] $D_n(z)$.

$$D_{-\nu}(z)D_{-\nu}(\pm z) = \frac{1}{\Gamma(\nu)} \sum_{n=0}^{+\infty} \frac{(\mp 1)^n D_n(z)^2}{n!(n+\nu)} \quad (61)$$

Formula (61) can also be found in Malyshev paper [7] or in its generalization version provided by Glasser, in [8], for unrelated arguments $y$ and $z$.

Furthermore, from corollaries 1 and 2, it is possible to extract other series expansions of the product of parabolic cylinder functions. If we let $\mu = 1$ in (43) and use the generating function of Laguerre polynomial $L_n^{-1}(z)$ in [5, p. 1002]

$$e^{-\frac{z^2}{2}\coth t} = e^{-\frac{z^2}{2}} \sum_{n=0}^{+\infty} L_n^{-1}(z^2) e^{-2nt}, \ t > 0 \quad (62)$$

we obtain, after integrating term by term and simplifying, the following series expansion of $D_{-\nu}(z)^2$

$$D_{-\nu}(z)^2 = \pi e^{-\frac{z^2}{2}} \sum_{n=0}^{+\infty} L_n^{-1}(z^2) \frac{2^{-\nu-n}(\nu)_n}{\Gamma(\frac{1+n+\nu}{2})^2} \quad (63)$$

Likewise, if we proceed as in (63), we can easily obtain

$$D_{-\nu}(z)D_{-\nu}(-z) = \pi e^{-\frac{z^2}{2}} \sum_{n=0}^{+\infty} (-1)^n L_n^{-1}(z^2) \frac{2^{-\nu-n}(\nu)_n}{\Gamma(\frac{1+n+\nu}{2})^2} \quad (64)$$

Identities (63) and (64) are valid for positive real number $z$ and $\nu \in \mathbb{C}$ or for all $z \in \mathbb{C}$ and $\nu$ is a negative integer.

## 4. Conclusions

In sum, we have proved that Laplace transform of the function $x \mapsto x^{\mu-1}e^{-\frac{1}{2}x^2}F(\nu,\mu;x^2)$, where $F = \Phi$ or $\Psi$, is explicit and equals the product of parabolic cylinder functions with different parameters $D_\nu(\pm z)D_{\nu+\mu-1}(z)$. These results lead also to Nicholson-type formulas and to new series expansions of the product $D_\nu(\pm z)D_{\nu+\mu-1}(z)$. In addition, we have established, in theorem 3, Laplace and Fourier transforms of the function $x \mapsto x^{2\nu-1}e^{-2x^2}F(\nu,\nu+\frac{1}{2};x^2)$, where again $F = \Phi$ or $\Psi$. Furthermore, the generalization of formulas in theorem 3 seems to do not involve confluent hypergeometric functions. I guess that their generalization would entail higher order confluent hypergeometric functions $_2F_2$.

Finally, the results of this paper can also be useful for the calculation of integrals involving the product of two, three or even four parabolic cylinder functions. References

---

[1] Hermite's polynomial $H_n$ is related to $D_n$ by $H_n(z) = 2^{\frac{n}{2}} e^{\frac{z^2}{2}} D_n(\sqrt{2}z)$.